\documentclass [11pt]{article}
\usepackage{amsfonts}
\usepackage{amsmath}
\usepackage{latexsym}
\usepackage{epsfig}
\usepackage{graphicx}
\def\e{\varepsilon}

\newtheorem{theorem}{Theorem}[section]
\newtheorem{lemma}[theorem]{Lemma}

\newtheorem{remark}{Remark}
\newtheorem{propos}[theorem]{Proposition}

\newtheorem{defin}{Definition}

\newenvironment{proof}{\noindent{\bf Proof.\ }}{\hfill$\Box$}
\date{}

\begin{document}
\title{Large  nearly regular induced subgraphs}
\author{
Noga Alon \thanks{ Schools of Mathematics and Computer Science,
Raymond and Beverly Sackler Faculty of Exact Sciences, Tel Aviv
University, Tel Aviv 69978, Israel and IAS, Princeton, NJ 08540,
USA. Email: {\tt nogaa@tau.ac.il}. Research supported in part by the
Israel Science Foundation, by a USA-Israeli BSF grant, 
by the Hermann Minkowski Minerva Center
for Geometry at Tel Aviv University and by the Von Neumann Fund.}
\and Michael Krivelevich \thanks{Department of Mathematics, Raymond
and Beverly Sackler Faculty of Exact Sciences, Tel Aviv University,
Tel Aviv 69978, Israel. E-mail: {\tt krivelev@post.tau.ac.il}. Research
supported in part by USA-Israel BSF Grant 2006322, and by grant
526/05 from the Israel Science Foundation.} \and  Benny Sudakov
\thanks{ Department of Mathematics,
UCLA,  Los Angeles, CA 90095 and Institute for Advanced Study, Princeton, 
NJ. Email: {\tt bsudakov@math.ucla.edu}.
Research supported in part by NSF CAREER award DMS-0546523, NSF
grants DMS-0355497 and DMS-0635607, by a USA-Israeli BSF grant, and by the 
State of New
Jersey.}    }

\maketitle
\begin{abstract}
For a real $c \geq 1$ and an integer $n$, let $f(n,c)$ denote the
maximum integer $f$ so that every graph on $n$ vertices contains 
an induced subgraph on at least $f$ vertices in which the 
maximum degree is at most $c$ times the minimum degree.
Thus, in particular, every graph on $n$ vertices contains
a regular induced subgraph on at least $f(n,1)$ vertices.
The problem of estimating $f(n,1)$ was posed long time ago
by Erd\H{o}s, Fajtlowicz and Staton. 
In this note we obtain the following upper and lower bounds for the 
asymptotic behavior of $f(n,c)$:

\noindent
(i) For fixed $c>2.1$, $n^{1-O(1/c)} \leq f(n,c)  \leq O(cn/\log n)$.

\noindent
(ii) For fixed $c=1+\varepsilon$ with $\varepsilon>0$ sufficiently small,
$f(n,c) \geq n^{\Omega(\varepsilon^2/ \ln (1/\varepsilon))}. $

\noindent
(iii) $\Omega (\ln n) \leq f(n,1) \leq O(n^{1/2} \ln^{3/4} n).$

\noindent
An analogous problem for not necessarily induced subgraphs is briefly considered
as well.
\end{abstract}

\section{Introduction}
All graphs considered here are finite and simple. For a graph
$G=(V,E)$, let $\Delta(G), \delta(G)$ and $d(G)=\frac{2|E|}{|V|}$
denote its maximum degree, minimum degree and average degree,
respectively. 
The {\em density} of $G$ is 
$p=|E|/{{|V|}\choose 2}$, clearly this is a number
between $0$ and $1$. For $U \subseteq V$, let $G[U]$ denote the subgraph of
$G$ induced on $U$.

\begin{defin}\label{def1}
A graph $G$ is {\em $c$-nearly regular} if
 $\Delta(G)\le c\cdot\delta(G)$.
\end{defin}

For a graph $G=(V,E)$ and a constant $c\ge 1$, let
$$
f(G,c)=\max\{|U|: G[U] \mbox{ is a $c$-nearly regular graph}\}\ .
$$
Define
$$
f(n,c)=\min\{f(G,c): |V(G)|=n\} .
$$
Thus, every graph $G$ on $n$ vertices contains a $c$-nearly regular
induced subgraph on at least $f(n,c)$ vertices. In particular, for $c=1$
every such $G$ contains a strictly regular induced subgraph on at
least $f(n,1)$ vertices.

The problem of estimating $f(n,1)$ was posed by 
Erd\H{o}s, Fajtlowicz and Staton (c.f. \cite{Er1}
or \cite{CG}, page 85). By the known estimates for Graph Ramsey numbers (c.f., e.g., \cite{GRS}), 
every graph on $n$ vertices contains either clique or independent set of order
$\Omega (\ln n)$. This implies that $f(n,1) \geq \Omega (\ln n)$. 
Erd\H{o}s, Fajtlowicz and Staton conjectured that the ratio $f(n,1) / \ln n$
tends to infinity as $n$ tends to infinity. We are unable to prove
or disprove this conjecture, and can only obtain several bounds,
listed in the following results. The first deals with
the case of large $c$.
\begin{propos}
\label{p11}
There exists an absolute constant $b$ so that for
$K\geq 2.1$, 
$f(n,K)\ge n^{1-b/K}$.
\end{propos}
The problem of obtaining a nontrivial lower bound for
values of $c$ close to $1$ is more interesting. Here we first 
deal with the case of graphs with positive density, and show that
any such graph must contain a nearly regular subgraph on a linear
number of vertices.
\begin{theorem}
\label{t12}
Let $\e>0$ be a small real, and  let $p$ satisfy $0 <p<1$.
Then, for every sufficiently large $n$, any graph $G=(V,E)$ on
$n$ vertices with density at least $p$ contains an induced 
$(1+\varepsilon)$-nearly regular subgraph on at least
$$
0.5 \left(\frac{\varepsilon}{6}\right)^{\frac{144}{\varepsilon^2} \ln (1/p)}\cdot  n
$$
vertices.
\end{theorem}
For general (possibly sparse) graphs we have the following:
\begin{theorem}
\label{t13}
Let $\e>0$ be a sufficiently small constant. Then
$$
f(n,1+\e) \ge n^{\frac{\e^2}{250\ln (1/\e)}}\,,
$$
for all sufficiently large $n$.
\end{theorem}
Our upper bounds for $f(n,c)$ are rather far from the lower bounds.
For the strictly regular case we prove the following.
\begin{theorem}
\label{t14}
$f(n,1) \leq O(n^{1/2}\log^{3/4}n)$.
\end{theorem}
This is a slight improvement of an earlier estimate
of Bollob\'as (c.f. \cite{CG}), who showed that for every $\epsilon>0$,
$f(n,1) \leq c(\epsilon)n^{1/2+\epsilon}$. 
For the nearly regular case we have:
\begin{propos}
\label{p15}
For every constant $K\ge 2$,
$f(n,K)\leq 7K\frac{n}{\log n}$.
\end{propos}

The lower bounds are proved in the next section, the upper bounds
are presented in Section 3. We conclude in Section 4 with
a few open problems and a brief discussion of an analogous problem for 
not necessarily induced  subgraphs.
Throughout this note we assume, whenever this is needed, that the 
number of vertices $n$ of the graphs discussed is sufficiently
large. To simplify the presentation, we make no attempt to optimize the absolute
constants, and omit all floor and ceiling
signs whenever these are not crucial. We also use the
following standard asymptotic notation:
for two functions $f(n)$, $g(n)$ of
a natural valued parameter $n$, we write $f(n)=o(g(n))$, whenever
$\lim_{n \rightarrow\infty} f(n)/g(n)=0$; $f(n)=O(g(n))$ if there
exists a constant $C>0$ such that $f(n)\le C g(n)$ for all $n$, and 
$f(n)=\Omega(g(n))$ if $g(n)=O(f(n))$.

\section{Lower bounds}
For a graph $G=(V,E)$ and a subset $U\subseteq V$, the number of
edges of $G$ spanned by $U$ in $G$ is denoted by $e_G(U)$; the
number of edges between disjoint subsets $U,W$ of vertices of $G$ is denoted
by $e_G(U,W)$.

\subsection{Large $c$}

In this subsection we prove Proposition \ref{p11}, that provides 
a lower bound for $f(n,c)$ when $c$ is
a relatively large constant.

We need the following rather standard
argument, allowing one to pass from a graph with a large average
degree to one with a large minimum degree.

\begin{propos}
\label{p21}
Let $K > 1,\alpha < 1/2$ be constants. Then every graph $G=(V,E)$
on $|V|=n$ vertices with $\Delta(G)\le Kd(G)$ contains an induced
$(K/\alpha)$-nearly regular subgraph $G^*$ with at least
$\frac{1-2\alpha}{K-2\alpha}n$ vertices and at least 
$\frac{K-2K\alpha}{2K-4\alpha} nd$ edges.
\end{propos}

\begin{proof}
Denote the average degree $d(G)$ of $G$ by $d$. 
We can obviously assume that $d>0$. Start with $G$
and delete repeatedly vertices of degree less than $\alpha d$ till
there are none left. Denote the resulting graph by $G^*$. Then
$G^*$ is an induced subgraph of $G$, satisfying $\Delta(G^*)\le
\Delta(G)$, $\delta(G^*)\ge \alpha d$, implying
$$
\Delta(G^*)\le \Delta(G)\le K d \le \frac{K}{\alpha} \delta(G^*)\,,
$$
and thus $G^*$ is a $(K/\alpha)$-nearly regular graph. We now
estimate the number of vertices of $G^*$. Denote the latter by $t$.
While creating $G^*$ from $G$, we deleted less than $(n-t)\alpha d$
edges, and thus
$$
\Delta(G^*)\ge d(G^*)> \frac{2(|E(G)|-(n-t)\alpha d)}{t}
=\frac{nd-2(n-t)\alpha d}{t}\ .
$$
But $\Delta(G^*)\le \Delta(G)\le Kd$, implying:
$$
\frac{nd-2(n-t)\alpha d}{t}\le Kd\ .
$$
Solving the above inequality for $t$, we get $t\ge
\frac{1-2\alpha}{K-2\alpha}n$, supplying the required lower bound for the number of
vertices of $G^*$. To bound the number of its edges  note that the number of vertices
deleted is $n-t$ and hence the number of edges deleted is at most $(n-t)\alpha d$,
leaving at least 
$$
\frac{1}{2}nd-(n-t)\alpha d=
\left(\frac{1}{2}-\alpha\right)nd+t\alpha d \geq
\left(\frac{1}{2}-\alpha\right)nd+\frac{1-2\alpha}{K-2\alpha}n \alpha d
=\frac{K-2K\alpha}{2K-4\alpha} nd,
$$
as needed.
\end{proof}

\begin{remark}\label{rem1}
It is instructive to observe that the above argument breaks down
completely for $\alpha\ge 1/2$. Therefore, when estimating $f(n,c)$ from
below for small $c$, in particular for $c<2$, we will adapt
a different strategy.
\end{remark}

We proceed with the following result, whose proof resembles
that of one of the  results in
\cite{ES}. 
\begin{propos}
\label{p22}
Let $K>1$ be a constant. Every graph $G=(V,E)$ on $|V|=|V(G)|=n$
vertices contains an induced subgraph $G^*$ on at least
$n^{1+\log_2\left(1-\frac{1}{K}\right)}$ vertices, for which
$\Delta(G^*)\le K d(G^*)$.
\end{propos}

For large $K$ the above estimate behaves like
$n^{1-\Theta\left(\frac{1}{K}\right)}$. Therefore, the assertions of
Proposition
\ref{p21}  and Proposition \ref{p22} imply that of Proposition 
\ref{p11}.

\noindent{\bf Proof of Proposition \ref{p21}.\ } Set $G_0=G$,
$k^*=\log_2n$. For $i=0,\dots,k^*$  repeat the following
loop. Set
$$
n_i=|V(G_i|,\quad \Delta_i=\Delta(G_i),\quad d_i=d(G_i)\ .
$$
If $\Delta_i\le Kd_i$, abort the loop. Otherwise delete
repeatedly vertices of degree at least $Kd_i/2$ from $G_i$ till
there are none left. Let $G_{i+1}$ be the resulting graph
and increment $i$.

Denote by $G^*$ the resulting graph of the above described process.
Observe that at iteration $i$ we delete at most
$|E(G_i)|/(Kd_i/2)=(n_id_i/2)(Kd_i/2)= n_i/K$ vertices, and thus
$n_{i+1}\ge (1-1/K)n_i$. It follows that
$$
|V(G^*)|\ge
\left(1-\frac{1}{K}\right)^{k^*}n=
n^{1+\log_2\left(1-\frac{1}{K}\right)}\ .
$$
If $G^*$ was created when the above loop was aborted due to
$\Delta(G_i)\le Kd(G_i)$, then obviously the obtained graph meets
the claim of the theorem. Otherwise, $G^*$ was obtained after $k^*$
iterations. At each such iteration $i$, we have $\Delta_{i+1}\le
Kd_i/2$ and $d_i \le \Delta_i/K$, implying $\Delta_{i+1}\le
\Delta_i/2$. Therefore, in this case
$$
\Delta(G^*)\le\Delta_0\,\cdot\, \left(\frac{1}{2}\right)^{k^*}<
n\cdot \left(\frac{1}{2}\right)^{\log_2n}=1\,,
$$
implying that $G^*$ has no edges and thus $\Delta(G^*)=d(G^*)=0$,
and $G^*$ can again serve as the required graph. 
\hspace*{\fill}\mbox{$\Box$}

\subsection{Small $c$}

Next we treat the more challenging case where the 
constant $c$ in $f(n,c)$ is very
close to 1. Throughout  this subsection $\varepsilon$ denotes a small 
positive real.

We start with several lemmas. 
\begin{lemma}
\label{l23}
Let $G=(V,E)$ be a graph on $n$ vertices 
with density $p$, and let $\varepsilon>0$. 
Then $G$ contains an induced subgraph $G'$ of density
$p' \geq p$, on a set  of 
$$
n' \geq \varepsilon^{\frac{2}{\varepsilon} \ln (1/p)} n
$$
vertices so that every set of $t \geq \varepsilon n'$ vertices
of $G'$ spans at most ${t \choose 2} p' (1+\e)$ edges.
\end{lemma}
\begin{proof}
Set $G_0=G$. For $i=0,1,\ldots$ repeat the
following loop. Set
$$
n_i=|V(G_i)|,\quad m_i=|e(G_i)|,\quad p_i=\frac{m_i}{{{n_i}\choose
2}}\ .
$$
If $G_i$ contains a subset $U_i\subseteq V(G_i)$ of at least
$\e n_i$ vertices such that
$e_{G_i}(U)\ge {{|U_i|}\choose 2}p_i(1+\e)$, then set
 $G_{i+1}:=G_i[U]$, $i:=i+1$.

Observe that after $k$ iterations of the above loop, the density $p_k$ of
the current graph $G_k$ satisfies: $p_k\ge (1+\e)^kp_0$. Thus, if the
loop is repeated at least $\frac{2}{\e} \ln (1/p)$ times, we have:
$$
p_k\ge (1+\e)^{\frac{2}{\e}\ln (1/p)}\,\cdot\, p
  > 1
 $$
-- a contradiction. It follows that the above process concludes after less
than $\frac{2}{\e} \ln (1/p)$ iterations. 
The resulting graph $G_k$ has $n_k$
vertices and $m_k$ edges. Observe that at each iteration the number
of vertices of the new graph is at least an $\e$-proportion of
the number of vertices of the previous graph. Therefore,
$$
n_k\ge \e^k |V(G)| \ge \e^{\frac{2}{\e} \ln (1/p)}n.
$$
We can thus take $G'=G_k$, $n'=n_k$ to complete the proof.
\end{proof}

\begin{lemma}
\label{l24}
Let $G=(V,E)$ be a graph on $|V|=n$ vertices with
$m=|E|$ edges and density $p=m/{{n}\choose 2}\ge n^{-a}$,
for some constant $0<a<1$. Suppose that 
\begin{equation}
\label{e21}
\mbox{ every } 
t\ge \e n 
\mbox{ vertices in } 
G 
\mbox{ span at most }
{t\choose 2}p(1+\e) 
\mbox{ edges. }
\end{equation}
Then for every subset $U\subseteq V$
of cardinality $|U|=\e n$ in $G$, there are at most
$$
\e n^2p(1+2\sqrt{\e})
$$
edges between $U$ and its complement in
$G$.
\end{lemma}

\begin{proof}
Assume that $U\subset V$ contradicts the above statement.
Denote:
$$
e_1=e_{G}(U), \quad e_2=e_{G}(U,V-U), \quad
e_3=e_{G}(V-U)\ .
$$
Then $e_1+e_2+e_3=m={{n}\choose 2}p$.

Choose uniformly at random a subset $X\subset V-U$ of cardinality
$|X|=x=\sqrt{\e}|V-U|=\sqrt{\e}(1-\e)n$. Then the expected
number of edges of $G$ spanned by $U\cup X$ is:
\begin{eqnarray*}
E[e_{G}(U\cup X)] &=& e_1+\frac{x}{|V-U|}e_2+
\frac{x(x-1)}{|V-U|(|V-U|-1)}e_3 \\
&>& e_1+\frac{x}{n-\e n}e_2+\frac{x^2}{(n-\e n)^2}
\left(1-\frac{1}{x}\right)e_3\\
&=& e_1+\sqrt{\e}e_2+\e e_3-O(n) =
e_1+\sqrt{\e}e_2+\e(m-e_1-e_2)- O(n)\\
&\ge& (\sqrt{\e}-\e)e_2+\e m-O(n)\ge
(\sqrt{\e}-\e)\e n^2p(1+2\sqrt{\e})+\e m-O(n)\\
&=&
(\sqrt{\e}-\e)\e(1+2\sqrt{\e})n^2p+\frac{\e n^2p}{2}-O(n)
=: A\ .
\end{eqnarray*}

On the other hand, by the assumption
on $G$, every such
set $U\cup X$ satisfies:
\begin{eqnarray*}
e_{G}(U\cup X) &\le & {{\e n+\sqrt{\e}(1-\e)n}\choose 2}
p(1+\e)\\
&\le& \frac{\e n^2p}{2}(1+\sqrt{\e}-\e)^2(1+\e)=
\frac{\e n^2p}{2}(1+2\sqrt{\e}-\e-2\e^{3/2}+\e^2)(1+\e)\\
&=& \frac{\e n^2p}{2}(1+2\sqrt{\e}+O(\e^{5/2})) =: B\ .
\end{eqnarray*}

Let us compare the asymptotic (in small $\e$) behavior of the two
quantities $A$ and $B$ defined above. We have:
\begin{eqnarray*}
A &=&
\frac{\e n^2p}{2}(1+(2\sqrt{\e}-2\e)(1+2\sqrt{\e}))-O(n)\\
&=& \frac{\e n^2p}{2}
(1+2\sqrt{\e}+2\e-4\e^{3/2})-O(n)\ .
\end{eqnarray*}
Since $n=o(n^2p)$, we have that $A>B$ for $\e$ small enough -- a contradiction.
\end{proof}

\begin{lemma}
\label{l25}
Let $G=(V,E)$ be a graph on $|V|=n$ vertices with
$m=|E|$ edges and density $p=m/{{n}\choose 2}\ge n^{-a}$,
for some constant $0<a<1$. Suppose that (\ref{e21}) holds.
Then $G$ contains an induced subgraph $G^*$
on at least $(1-\e-2 \sqrt {\e})n>n/2$ vertices with maximum
degree $\Delta(G^*) \leq (1+3\sqrt{\e}) pn$ and minimum degree
$\delta(G^*) \geq (1-2 \sqrt {\e}) np$.
In particular, $G^*$ is $c$-nearly regular for $c=(1+6 \sqrt
{\e})$.
\end{lemma}

\begin{proof}
Let
$U$ be a set of $\e n$ vertices of highest degrees in $G$ (ties
are broken arbitrarily). Set $H=G[V-U]$. We claim that all
vertex degrees in $H$ are at most $np(1+3\sqrt{\e})$. If
this is not so, then the degrees of all vertices of $U$ in $G$ are
at least $n p(1+3\sqrt{\e})$, implying 
(through condition (\ref{e21})):
\begin{eqnarray*}
e_{G}(U,V-U)&\ge& \e n^2p(1+3\sqrt{\e})-2e_{G}(U) \ge
\e n^2 p(1+3\sqrt{\e})-\e^2 n^2 p(1+\e)\\
&>& \e n^2 p(1+2\sqrt{\e})\,,
\end{eqnarray*}
for small enough $\e$, thus contradicting Lemma  \ref{l24}.
Therefore, $H$ is an
induced subgraph of $G$  on $|V(H)|=(1-\e)n$ vertices, of
maximum degree $\Delta(H)\le n p(1+3\sqrt{\e})$, still
satisfying condition (\ref{e21}), and having
\begin{eqnarray}
\label{eq1}
|E(H)|&\ge&
{{n}\choose 2}p-\frac{\e^2 n^2}{2}p(1+\e)-\e n^2p(1+2\sqrt{\e})\nonumber\\
&\ge& {{n}\choose 2}p-2\e n^2p
\end{eqnarray}
edges.

We now delete from $H$ repeatedly vertices of degree less than
$n p(1-2 \sqrt{\e})$, until there are no such vertices, or until
we have deleted $2\sqrt{\e} n$ of them. Assume the latter case
happens, and denote the set of $2\sqrt{\e}n$ deleted vertices by
$W$. Then the set $V(H)-W$ has $|V(H)-W|=(1-\e-2\sqrt{\e})n$
vertices and by (\ref{eq1}) spans at least

\begin{eqnarray*}
e_H(V(H)-W) &\geq& |E(H)|-2\sqrt{\e}n\,\cdot\, n p(1-2\sqrt{\e}) \ge
{{n}\choose 2}p-2\e n^2p-2\sqrt{\e}n^2 p(1-2\sqrt{\e})\\
&=& {{n}\choose 2} p-2\sqrt{\e}n^2 p+2\e n^2p=
\frac{n^2p}{2}\left(1-4\sqrt{\e}+4\e\right)-O(np)
\end{eqnarray*}
edges. On the other hand, by condition (\ref{e21}), the set $V(H)-W$
satisfies:
\begin{eqnarray*}
e_H(V(H)-W) &\le& {{(1-2\sqrt{\e}-\e)n}\choose 2}p(1+\e)
\le (1-2\sqrt{\e}-\e)^2(1+\e)\frac{n^2p}{2}\\
&=&\frac{n^2p}{2}\left(1-4\sqrt{\e}+3\e+O(\e^{3/2})\right).
\end{eqnarray*}
Comparing the above two estimates for $e_H(V(H)-W)$ we get a
contradiction for small enough $\e$.

It follows that the above deletion process stops before
$2\sqrt{\e}n$ vertices have been deleted. Denote the resulting
graph by $G^*$. Then $G^*$ has $|V(G^*)|\ge
n-\e n-2\sqrt{\e}n > \frac{n}{2}$ vertices, has
maximum degree $\Delta(G^*)\le np(1+3\sqrt{\e})$ and minimum
degree $\delta(G^*)\ge n p(1-2 \sqrt{\e})$. Hence
$$
\Delta(G^*)\le \frac{1+3\sqrt{\e}}{1-2 \sqrt{\e}}\delta(G^*)<
(1+6\sqrt{\e})\delta(G^*),
$$
completing the proof of the lemma.
\end{proof}

We are now ready to prove Theorems \ref{t12} and \ref{t13}.
\medskip

\noindent
{\bf Proof of Theorem \ref{t12}.}\, 
By Lemma \ref{l23} (with $\frac{\e^2}{36}$ playing the role of $\e$),
$G$ contains an induced subgraph $G'$ of density $p' \geq p$ on
$$
n' \geq \left(\frac{\e^2}{36}\right)^{\frac{72}{\e^2} \ln (1/p)} \cdot n
=\left(\frac{\e}{6}\right)^{\frac{144}{\e^2} \ln (1/p)} \cdot n
$$
vertices, such that every set of $t \geq \frac{\varepsilon^2}{36} n'$ 
vertices of $G'$ spans at most
${t \choose 2}p'(1+\frac{\e^2}{36})$ edges. Since $p' \geq p > n^{-1/2}$ (as $p>0$ is a
constant and $n$ is large),
Lemma \ref{l25} implies that $G'$ contains an induced subgraph $G^{*}$ on at least $n'/2$
vertices which is $c$-nearly-regular for $c=1+6 \sqrt { \frac{\e^2}{36}} =1+\e$,
as needed. \hfill $\Box$
\newpage

\noindent
{\bf Proof of Theorem \ref{t13}.}\, 

Set 
$$
\e_0=\frac{\varepsilon^2}{36},\quad a=\frac{\e_0}{3\ln(1/\e_0)}\ .
$$
Let $G=(V,E)$ be a graph on $n$ vertices. Denote by
$p=|E|/{n\choose 2}$ the density of $G$.
The average degree of $G$ is at most $np$ and by Tur\'an's theorem
$G$ contains an independent set $U$ of
size $n/(np+1)$. Therefore we can assume that
$p\ge n^{-a}$, as otherwise 
\begin{eqnarray*}
\frac{n}{np+1}&\ge& \frac{n}{n^{1-a}+1}\ge (1-o(1))n^{a} =
(1-o(1))n^{\frac{\e_0}{3\ln(1/\e_0)}}\\
&=& (1-o(1))n^{\frac{\e^2}{108\ln(36/\e^2)}} >
n^{\frac{\e^2}{250\ln(1/\e)}}\,,
\end{eqnarray*}
where here we used the assumption that $\e$ is sufficiently small and
$n$ is sufficiently large. This gives an induced $0$-regular subgraph of $G$, and we can thus indeed
assume that $p\ge n^{-a}$.

By Lemma \ref{l23} $G$ contains an induced subgraph $G'$ of density $p' \geq p$
on 
$$
n' \geq \e_0^{\frac{2}{\e_0} \ln(n^a)} \cdot n = \e_0^{\frac{2}{3 \ln (1/\e_0)} \ln n} \cdot n=
n^{1/3}
$$ 
vertices, in which the density of the induced subgraph on any
set of at least $\e_0n'$ vertices does not exceed
$p' (1+\e_0)$. By Lemma \ref{l25}, $G'$ (and hence $G$) contains an induced
subgraph on at least $0.5 n' \geq 0.5 n^{1/3}$ vertices, which is
$(1+6 \sqrt {\e_0})=(1+\e)$-nearly regular, completing the proof. 
\hfill $\Box$

\begin{remark}
Note that we have actually proved the following result, which is
stronger than the assertion of Theorem \ref{t13}: Every graph on $n$ vertices
contains either an independent set of size at least $n^{\frac{\e^2}{250 \ln (1/\e)}}$,
or a $(1+\e)$-nearly regular induced subgraph on at least $0.5 n^{1/3}$ vertices.
\end{remark}

\section{Upper bounds}

\subsection{The strictly regular case}

{\bf Proof of Theorem \ref{t14}.}\,
Fix an integer $k$ satisfying
$$
k \geq Cn^{1/2}\ln^{3/4}n \,
$$
where $C>0$ is a sufficiently large constant to be set later.

We will work with the following model of random graphs on $n$
vertices which we denote by $G(n,\bar{p})$. Let
$\bar{p}=(p_1,\ldots,p_n)$, where
$$
p_i=\frac{1}{4}+\frac{i}{2n},\quad i=1,\ldots,n\ .
$$
Then $G(n,\bar{p})$ is the probability space of graphs with vertex
set $[n]=\{1,\ldots,n\}$, where for each pair $1\le i\ne j\le n$,
$(i,j)$ is an edge of $G(n,\bar{p})$ with probability $p_ip_j$,
independently of all other pairs. Notice that the probability of
each individual pair $(i,j)$ to be an edge of $G(n,\bar{p})$ is
strictly between $1/16$ and $9/16$.

\begin{propos}\label{prop1}
Let $X_1,\ldots, X_t$ be independent Bernoulli random variables,
where $Pr[X_i=1]=\rho_i,\ i=1,\ldots,t$. Let $X=X_1+\ldots+X_t$.
Assume that $1/16\le \rho_i\le 9/16$ for $i=1,\ldots,t$. Then for
every integer $0\le s\le t$, $Pr[X=s]\le c_0/\sqrt{t}$, for some
absolute constant $c_0>0$.
\end{propos}

\begin{proof}
For every $1\le i\le t$, we represent $X_i$ as a product $X_i=Y_i\cdot Z_i$, where
$\{Y_i,Z_i\}$ is a collection of mutually independent Bernoulli
random variables defined by $Pr[Y_i=1]=9/16$,
$Pr[Z_i=1]=\frac{16}{9}\rho_i$, $i=1,\ldots, t$.

Set $I_0=\{1\le i\le t: Z_i=1\}$. Since $\rho_i\ge 1/16$ for
all $1\le i\le t$, we have that  $Pr[Z_i=1] \geq 1/9$ and $E[|I_0|] \geq t/9$.   
By standard large deviation arguments,
$|I_0|\ge t/10$ with probability $1-o(1/\sqrt{t})$. Thus
\begin{eqnarray*}
Pr[X=s] &=& \sum_{I\subseteq [t]} Pr[I_0=I]\cdot Pr[\sum_{i\in
I}Y_i=s]= \sum_{I\subseteq [t]}Pr [I_0=I] Pr[B(|I|, 9/16)=s]\\
&\le& Pr[|I_0|< t/10] + \sum_{|I|\ge t/10} Pr[I_0=I] Pr[B(|I|,9/16)=s]
\ .
\end{eqnarray*}
Here $B(n,p)$ denotes the binomial random variable with parameters
$n$ and $p$. From known estimates on binomial random variables, we
obtain that $Pr[B(r,9/16)=s]\le \frac{c}{\sqrt{r}}\leq \frac{4c}{\sqrt{t}}$ for every $r\ge
t/10$, where $c>0$ is an absolute constant. Plugging this estimate
into the inequality above, we get the claimed result.
\end{proof}

\begin{lemma}\label{lem1}
Let $U\subseteq [n]$ be a fixed set of $|U|=k$ vertices. Then the
probability that in $G(n,\bar{p})$, with $\bar{p}$ as defined above, 
the graph $G[U]$ is a regular graph is at most
$n\left(\frac{c_1}{k}\right)^{k/2}$ for some absolute constant
$c_1>0$.
\end{lemma}

\begin{proof}
Fix the degree of regularity $d$ of the regular subgraph $G[U]$
(this can be done in $n$ ways).  Let $U=\{u_1, u_2 \ldots ,u_k\}$.
We bound the probability that the induced subgraph $G[U]$
is $d$-regular as follows.

Expose the edges of $G[U]$ by first exposing the edges from $u_1$
to $U-\{u_1\}$, then from $u_2$ to $U-\{u_1,u_2\}$, etc. If
vertex $u_i$ gets $t_i$ neighbors in $\{u_1,\ldots,u_{i-1}\}$, then
$u_i$ should have exactly $d-t_i$ neighbors in $\{u_{i+1},\ldots,
u_{k_0}\}$. Recall that all edge probabilities in $G(n,\bar{p})$ are
between $1/16$ and $9/16$. Thus Proposition \ref{prop1} applies, and
the probability of the latter event is at most $c_0/\sqrt{k-i}$.
Multiplying these probabilities for $i=1,\ldots,k$, we derive that
the probability that $U$ is a $d$-regular graph is at
most
$$
\prod_{i=1}^{k-1} \frac{c_0}{\sqrt{k-i}} =
\frac{(c_0)^{k-1}}{\sqrt{(k-1)!}} \le
\left(\frac{c_1}{k}\right)^{k/2}
$$
where $c_1>0$ is an absolute constant and 
the last inequality follows by applying the Stirling formula.
\end{proof}

\medskip

We can now complete the proof of Theorem \ref{t14}. 
We bound the probability 
that in $G(n,\bar{p})$ there
exists a set $U$ of size $|U|=k$, such that
$U=(u_1<\ldots<u_{k/\ln k}=a<\ldots < u_{k-k/\ln k}=b<\ldots<u_k)$ 
and $G[U]$
is a regular graph by considering two possible cases depending
on the difference $t=b-a$ between $a$ and $b$.

\noindent
{\bf Case 1:} $t \leq c_2k^{3/2}$, where $c_2>0$ is a small positive 
constant to be determined later.

The probability that there exists such a $U$ is at most:
\begin{equation}
\label{e99}
{n \choose {2k/\ln k}} n^2  {t\choose{k-2k/\ln k}}
\cdot n \left(\frac{c_1}{k}\right)^{k/2}.
\end{equation}
Indeed, there are less than  ${n \choose {2k/\ln k}} $ ways 
to choose the vertices 
$$u_1, \ldots ,u_{k/\ln k-1}, u_{k-k/\ln k+1}, \ldots, u_k,$$
less than $n^2$ ways to choose the vertices $a$ and $b$
so that the difference between them, $t$, is at most
$c_2 k^{3/2}$, and  less than
${t\choose{k-2k/\ln k}}$ ways to choose the vertices 
$u_{k/\ln k+1}, \ldots ,u_{k-k/\ln k-1}$. For each such choice,
the probability that the induced subgraph on $\{u_1, u_2 ,
\ldots ,u_k\}$ is regular is at most 
$n \left(\frac{c_1}{k}\right)^{k/2}$, 
by Lemma \ref{lem1}.

A simple computation shows that the expression in (\ref{e99}) is
(much) smaller than, say, $1/n^2$ for an appropriate choice of
$c_2$. Indeed, since $n \leq k^2$ and 
${t\choose{k-2k/\ln k}} \leq {t \choose k} \leq (et/k)^k$, this expression
is at most
$$
n^3 k^{4k/\ln k} \left(\frac{et}{k}\right)^k
\left(\frac{c_1}{k}\right)^{k/2}
\leq n^3 e^{5k} c_2^k c_1^{k/2}=n^3 [e^{10} c_2^2 c_1]^{k/2},
$$
implying the required estimate by choosing, for example,
$c_2 = \frac{1}{2 e^5 \sqrt {c_1}}$.

\noindent
{\bf Case 2:}\, $t \geq c_2 k^{3/2}$.

Let $U=(u_1<\ldots<u_{k/\ln k}=a<\ldots <u_{k-k/\ln k}=b<\ldots 
< u_k)$.
Denote the first (smallest) $k/2$ vertices of $U$ by $U_1$, and the
last $k/2$ vertices by $U_2$. Observe that if $G[U]$ is a regular
graph, then the two induced subgraphs $G[U_1]$ and $G[U_2]$ have the
same average degree. This is highly improbable. Indeed, by the
definition of $G(n,\bar{p})$ the probability of each pair inside
$U_1$ to become an edge is strictly less than the probability of
each pair inside $U_2$ to become an edge. In addition, since
$b-a=t$, the probability of each pair $i,j\in U_1$ with $i <j, 
i \le a$
to be an edge is less than the probability of each pair $i',j'\in
U_2$ with $i' < j', j' \ge b$ to be an edge by 
$\Omega(t/n)$ -- this is
because in this case
$$
\left(\frac{1}{4}+\frac{i'}{2n}\right)\left(\frac{1}{4}+\frac{j'}{2n}\right)
-
\left(\frac{1}{4}+\frac{i}{2n}\right)\left(\frac{1}{4}+\frac{j}{2n}\right)
\geq \Omega\left(\frac{i'+j'-i-j}{2n}\right) 
\geq \Omega\left(\frac{j'-i}{2n}\right) \geq
\Omega \left(\frac{t}{n}\right).
$$
It follows that the expected number of edges inside $U_2$ exceeds
that inside $U_1$ by $\Omega\left(\frac{k^2t}{n \ln k}\right)$. By
Chernoff's Inequality (c.f., e.g., \cite{AS}, Appendix A), 
we obtain that the probability that
$G[U_1]$ and $G[U_2]$ have the same average degree in $G(n,\bar{p})$
is at most $\exp\{-c_3k^2t^2/(n^2 \ln^2 k)\}$ 
for some absolute constant
$c_3>0$.

Thus, recalling our assumption on $t$ we derive that the probability
that there is $U$ as above for which $G[U]$ is a regular graph is at
most
$$
{n\choose k}e^{-\frac{c_3k^2t^2}{n^2 \ln^2 k}}\le
\left(\frac{en}{k}\right)^k e^{-\frac{c_3c_2^2k^5}{n^2 \ln^2 k}} =
\left(\frac{en}{k}e^{-\frac{c_3c_2^2k^4}{n^2 \ln^2 k}}\right)^k\ .
$$
Finally, as $k \geq Cn^{1/2}\ln^{3/4}n$ and $\ln k=\Theta( \ln n)$ 
we can choose $C>0$ to be large
enough so that the above expression is (much) smaller than $1/n^2$. 

Combining the two cases we conclude that the probability that for 
any fixed $k$ which is at least $Cn^{1/2} \ln^{3/4} n$ our graph 
contains an induced regular subgraph on $k$ vertices is smaller than
$2/n^2<1/n$, and as there are less than $n$ choices for $k$, this
shows that with positive probability the graph contains no such
subgraph.
This completes
the proof of Theorem \ref{t14}.  \hfill $\Box$
\bigskip

\subsection{The nearly regular case}
\noindent
{\bf Proof of Proposition \ref{p15}.}\,
We will prove that for every (large) $n$ there exists a graph $G$ on $n$
vertices in which every $K$-nearly regular induced subgraph has at
most $7Kn/\log n$ vertices.

 Assume first that $n$ is of the form $n=(s+1)2^s$ for a positive
 integer $s$. Notice that $s=(1-o(1))\log_2n$. Take a set $V$ of $n$
 vertices and partition it into $s+1$ disjoint equally sized subsets
 $V_0,\ldots, V_{s}$, $|V_i|=2^s$. Now we define $G$ as follows. For
 $i=0,\ldots,s$ the set $V_i$ spans $2^{s-i}$ disjoint cliques of size
 $2^i$ each. There are no edges between the cliques inside $V_i$ and
 no edges between distinct subsets $V_i\ne V_j$ in $G$.

 Assume now that a subset $U\subseteq V$ spans a graph $G[U]$
 satisfying: $\delta(G[U])=d$, $\Delta(G[U])\le Kd$. Observe that
 the degrees of all vertices from $V_i$ in $G$ are $2^i-1$, and thus if
 $2^i-1<d$, then $U\cap V_i=\emptyset$. Let now $2^i\ge d+1$. Since 
$\Delta(G[U])\le Kd$, 
 $U$ has at most $Kd+1$ vertices in each clique spanned by $V_i$,
 implying $|U\cap V_i|\le 2^{s-i}(Kd+1)$. Therefore:
\begin{eqnarray*}
|U| &=& \sum_{i=0}^s |U\cap V_i|=\sum_{i: 2^i\ge d+1} |U\cap V_i|
\le\sum_{i=\lceil \log_2(d+1)\rceil}^s2^{s-i}(Kd+1)
< (Kd+1)2^{s-\lceil \log_2(d+1)\rceil+1}\\
&\le&
\frac{2(Kd+1)}{d+1}2^s=\frac{2(Kd+1)}{d+1}\,\frac{n}{s+1} \leq
(2+o(1)) K\,\frac{n}{\log n}\,,
\end{eqnarray*}
implying the desired result.

For $n$ not of the form $n=(s+1)2^s$, choose a minimal $s$
satisfying $(s+1)2^s\ge n$. Let $n'=(s+1)2^s$. It is easy to verify
that $n'\leq 3n$. Now we can apply the above construction to
create a graph $G'$ on $n'$ vertices in which every $K$-nearly
induced subgraph has at most $(2+o(1)) Kn'/\log n' \leq 7Kn/\log n$ vertices, and then
take $G$ to be an arbitrary induced subgraph of $G'$ on exactly $n$
vertices.  \hfill $\Box$

\section{Open Problems}

The most intriguing open problem is that of obtaining a better estimate
for $f(n,1)$. In particular, the conjecture of Erd\H{o}s, Fajtlowicz and Staton that
$f(n,1)/\ln n$ tends to infinity as $n$ tends to infinity remains open. The values of
$f(n,1)$ for $n \leq 17$ have been determined by the authors of \cite{FMRS} and by
McKay, and these are indeed larger than the bounds that follow from the corresponding 
Ramsey numbers.

Our upper and lower bounds for $f(n,c)$ for $c>1$ are also rather far from each other,
and it will be nice to understand the behavior of this function better.

One can also study a
variant of the problems considered here that deals with not necessarily induced
subgraphs. Of course, every graph contains a regular subgraph on all vertices
(the subgraph with no edges), and hence in this case it is natural to look for 
regular or nearly regular subgraphs with a large number of {\em edges}.
For every two positive integers $n,m$ with 
$m \leq {n \choose 2}$ and a real $c \geq
1$, let $g(n,m,c)$
denote the largest $g$ so that every graph with $n$ vertices and $m$ edges contains
a (not necessarily induced) $c$-nearly regular subgraph with at least $g$ {\em edges}.  
The problems of determining or estimating
the behavior of this function seems interesting.  Here we can establish
tighter estimates  than the ones obtained for the induced case.

Consider first the case $c=1$. 
Since the complete graph on $n$
vertices can be covered by $n$ matchings (and by $(n-1)$ for even $n$), it follows
that $g(n,m,1) \geq m/n$, since every graph with $m$ edges contains a matching
of size at least $m/n$. The star $K_{1,n-1}$ shows that for some values of $m$ and $n$
this is essentially tight, and that $g(n,m,1)=1$ for all $1 \leq m<n$.
By a simple application of 
Szemer\'edi's Regularity Lemma it can be shown (see \cite{PRS}) that for every fixed $p>0$ there is
a $\delta=\delta(p)>0$ so that $g(n,pn^2,1) \geq \delta n^2$. 
This bound was significantly improved by 
R\"odl and Wysocka \cite{RW}, who proved that every graph with $n$ vertices and 
$pn^2$ edges contains an $r$-regular subgraph with $r\geq \alpha p^3n$ for some
positive constant $\alpha$.

For a larger constant $c$ observe, first, that
complete bipartite graphs 
show that for $m \geq n$, $g(n,m,c) \leq O(c(m/n)^2)=O(c d^2)$,
where $d=2m/n$ is the average degree of a graph with $n$ vertices
and $m$ edges.  Indeed, a complete bipartite
graph $K_{k,n-k}$ with $k \leq n/2$ 
has average degree  $d=\Theta(k)$.  Every
$c$-nearly-regular subgraph in it has minimum
degree at most $k$, and hence maximum degree at most $ck$. Thus it cannot have more
than $k \cdot ck=ck^2$ edges.  Therefore, for every fixed $c>1$ 
there exists some
$C=C(c)$ so that $g(m,n,c) \leq C(m/n)^2$ for all $m>n$. We can show that for $c>2$ 
this is tight, up to a
constant factor; namely, for any $c>2$ there is  a $b=b(c) >0$ so that
$g(n,m,c) \geq b (m/n)^2$ for all $m>n$. For simplicity we present the proof only for
$c=5$, the proof for any other $c>2$ is similar.
\begin{theorem}
\label{t41}
Let $G=(V,E)$ be a graph with $|V|=n$ vertices, $|E|=m>n$ edges and average degree
$d=d(G)=2m/n$. Then $G$ contains a $5$-nearly regular subgraph with at least
$\frac{d^2}{2^{12}}$ edges.
\end{theorem}
\begin{proof}
We apply the method of Pyber in \cite{Py}, 
together with a few extra twists.
Clearly we may assume that $d \geq 2^6$.
First omit from $G$ repeatedly vertices of degree smaller than $d/2$, as long as there
are such vertices. As this process can only increase the average degree, it ends with a
nonempty graph $G'$ with minimum degree at least $d/2$. Now take a spanning bipartite
subgraph of $G'$ with the maximum number of edges. It is easy and well known that the
degree of every vertex in this bipartite subgraph is at least half its degree in $G'$,
giving a bipartite graph $H$ with minimum degree at least $d/4$. 
Put $H_1=H$.
Let $A$ and $B$ denote the two vertex classes of $H$, where $|A| \geq |B|$. 
Let $A_1 \subseteq A$ be a nonempty subset of $A$ which satisfies
$|N_{H_1}(A_1)| \leq |A_1|$ and $A_1$ is minimal with respect to containment
(subject to the condition above and to being nonempty).  Clearly there is such
an $A_1$, as $|N_{H_1}(A)| \leq |A|$ and $|N_{H_1}(v)| \geq d/4>1$ for all $v \in A$.
Put $N_{H_1}(A_1)=B_1$ and note that by the minimality of $A_1$, $|A_1|=|B_1|$.
By the minimality, again, and by Hall's theorem,
there is a matching $M_1$ saturating $A_1$ and $B_1$. Let $H_2$ be the
graph obtained from $H_1$ by deleting all edges of $M_1$.
Now let $A_2 \subseteq A_1$ be a nonempty, minimal subset of
$A_1$ satisfying $|N_{H_2}(A_2)| \leq |A_2|$. As before, there is such a set,
as $|N_{H_2}(A_1)| \leq |N_{H_1}(A_1)|=|A_1|$. The minimality shows, again, that
in fact $|N_{H_2}(A_2)| = |A_2|$, and that there is a matching $M_2$ saturating
$A_2$ and $N_{H_2}(A_2)=B_2$.  Proceeding in this manner we define a sequence
of sets 
$$
\emptyset\ne A_{d/4} \subseteq A_{d/4-1} \subseteq \cdots \subseteq A_2 \subseteq A_1
\subseteq A
$$
and
$$
\emptyset\ne B_{d/4} \subseteq B_{d/4-1} \subseteq \cdots \subseteq B_2 \subseteq B_1
\subseteq B
$$
where $|A_i|=|B_i|$ for all $i$, and a sequence of pairwise edge-disjoint matchings
$M_{d/4}, \cdots , M_2, M_1$, where $M_i$ is a perfect matching between $A_i$ and $B_i$.
Note that indeed this process does not terminate before these $d/4$ phases, as
initially all degrees in $H$ are at least $d/4$, and with the omission of
each matching the degrees drop by $1$.

For convenience we assume, from now on, that $d$ is a power of $2$ (otherwise, simply
consider only the first $d'/4$ sets $A_i,B_i$ and matchings $M_i$, where
$d'>d/2$ is the largest power of $2$ that does not exceed $d$).
Note that $|A_{d/8}| >d/8$, since every vertex of $B_{d/4}$ is incident with an edge
of each of the matchings $M_i$ for $d/8 \leq i \leq d/4$, and all these edges are
incident with vertices of $A_{d/8}$.
We consider two possible cases.

\noindent
{\bf Case 1:}\, For every $i$, $ 0 \leq i \leq \log_2 d -4$, 
$$
|A_{d/2^{i+4}}| > 2 |A_{d/2^{i+3}}|.
$$
In this case, 
$$
|A_1| > 2^{\log_2 d-3} |A_{d/8}| \geq \frac{d^2}{64},
$$
and the matching  $M_1$ is a regular subgraph with more than $d^2/64$ edges, supplying
the desired result (with room to spare).

\noindent
{\bf Case 2:}\,  There is an $i$, $0 \leq i \leq \log_2 d -4$, such that
$$
|A_{d/2^{i+4}}| \leq 2 |A_{d/2^{i+3}}|.
$$
In this case, take the minimum $i$ for which this holds. Then
$$
|A_{d/2^{i+3}}| > 2^i |A_{d/8}| \geq 2^{i-3} d.
$$
Let $H'$ be the graph consisting of the $\frac{d}{2^{i+4}}$ matchings
$M_j$ for $\frac{d}{2^{i+4}} \leq j <\frac{d}{2^{i+3}}$. The vertices of $H'$
are all those saturated by the largest matching among those, 
namely $M_{\frac{d}{2^{i+4}}}$.
Then the maximum degree in
$H'$ is exactly $\frac{d}{2^{i+4}}$ (as every vertex of $A_{d/2^{i+3}}$ has that
degree), and the average degree is at least half of that, since each of the 
$\frac{d}{2^{i+4}}$ matchings $M_j$ above is of size at least 
half that
of the largest one, which is spanning. 
As in $H'$ the degree of every vertex of $A_{d/2^{i+3}}$ is exactly
$\frac{d}{2^{i+4}}$, the total number of edges of $H'$ is at least
$$ 
|A_{d/2^{i+3}}| \cdot \frac{d}{2^{i+4}} 
\geq 2^{i-3} d \cdot \frac{d}{2^{i+4}}=2^{-7} d^2.
$$

Thus, $H'$ is a graph with maximum degree that exceeds the average degree
by a factor of at most $K=2$. We can now apply Proposition \ref{p21} with
$K=2$ and $\alpha=0.4$ to conclude that $H'$ contains a $K/\alpha=5$-nearly
regular subgraph with at least 
$$
\frac{K-2K\alpha}{2K-4\alpha} |E(H')| \geq \frac{1}{6}2^{-7} d^2
$$
edges, completing the proof. 
\end{proof}

\noindent
{\bf Acknowledgment.} \,
We thank Domingos Dellamonica and Vojta R\"odl for pointing an error 
in an early version of the paper.

\end{document}